\documentclass[12pt]{article}

\usepackage{amsfonts}

\newtheorem{thm}{Theorem}[section]
\newtheorem{lemma}[thm]{Lemma}
\newtheorem{cor}[thm]{Corollary}

\newcommand{\proof
}{\par\medskip\noindent {\bf Proof.\ \ }}

\newcommand{\be}{\begin{equation}}
\newcommand{\ee}{\end{equation}}
\newcommand{\openbox}{\leavevmode
  \hbox to8pt{\hfil\vrule\vbox to6pt{\hrule width6pt\vfil\hrule}\vrule}}

\newcommand{\qed}{\hbox to5pt{ } \hfill \openbox\bigskip\medskip}

\newcommand{\cF}{\mbox{$\cal F$}}

\newcommand{\N}{\mathbb N}
\newcommand{\Z}{\mathbb Z}
\newcommand{\Q}{\mathbb Q}

\title{Linear equations for the number of intervals which are isomorphic 
with Boolean lattices and the Dehn--Sommerville equations}
\author{G\'abor Heged\"{u}s
\\{\normalsize Johann Radon Institute for Computational and Applied Mathematics}
}

\begin{document}
\footnotetext{Research supported in part by OTKA grant K77476}


\maketitle

\begin{abstract}
Let $P$ be a finite poset. Let $L:=J(P)$ denote the lattice of order ideals of $P$. Let $b_i(L)$ denote the number of Boolean intervals of $L$ of rank $i$.

We construct a simple graph $G(P)$ from our poset $P$. Denote by $f_i(P)$ the number of the cliques $K_{i+1}$, contained in the graph $G(P)$.

Our main results are some linear equations connecting the numbers $f_i(P)$ and $b_i(L)$.

We reprove the Dehn--Sommerville equations for simplicial polytopes. 

In our proof we use free resolutions and the theory of Stanley--Reisner rings.
\end{abstract}
\medskip

\section{Introduction}

In my article I use the following proof technique:

Let $R$ denote the polynomial ring $\Q[x_1,\ldots ,x_n]$. Let $I$ be a 
monomial ideal. If we know a graded free resolution of the quotient module 
$M:=R/I$: 
$$
0\longrightarrow F_n \longrightarrow \ldots \longrightarrow F_1 \longrightarrow M \longrightarrow 0
$$
then we can compute the Hilbert function (or the  Hilbert polynomial) of $M=R/I$. 

On the other hand, we can compute the Hilbert function of the module $M=R/I$ by counting standard monomials.

If we compare these computations we get a new equation,
and we can prove some linear equations for the combinatorial invariants of the monomial ideal $I$.
  
We give here two applications in lattice theory and in the theory of polytopes.
 
Throughout this paper we use the following notations.

Let $(P,\preceq_P)$ be a fixed poset. We say that $J\subseteq  P$ is an {\em order ideal} of $P$, if $a\in J$ and $b\preceq_P a$, then $b\in J$.
Let $L:=J(P)$ denote the distributive lattice of order ideals of $P$. 

We say that a lattice $B$ is a {\em Boolean lattice}, if $B$ is distributive, $B$ has $0$ and $1$ and each $a\in B$ has a complement $a'\in B$.

Let $M$ denote an arbitrary finite distributive lattice. Let $l,k\in M$ with $l\leq k$. Then the set 
$$
[l,k]:=\{m\in M:~ l\leq m\leq k\}\subseteq M
$$
is called an {\em interval} in $M$.
Let $b_i(M)$ denote the number of intervals of $M$, which are isomorphic to the Boolean lattice of rank $i$.

In our main result we describe some linear equations for the numbers $b_i(L)$.

We can state our results in a more compact form if we associate the following graphs $G(P)$ to the poset $P$.

Let $P=\{q_1,\ldots ,q_p\}$ be a finite poset, $|P|=p$. We define a simple graph $G(P)$ as follows: let the vertex set of $G(P)$ be the disjoint union $\{x_1,\ldots ,x_p\}\cup \{y_1,\ldots ,y_p\}$. Define the edge set of $G(P)$ as
\begin{equation} \label{edge}
 \{\{x_i,x_j\}:~ 1 \leq i<j\leq p\}\cup  \{\{y_i,y_j\}:~ 1 \leq i<j\leq p\}\cup \{\{x_i,y_j\}:~ q_i\not\leq q_j\}.
\end{equation}

Let $j$ be a nonnegative integer. Denote by $f_j(P)$ the number of the cliques $K_{j+1}$, contained in the graph $G(P)$.
 
In particular, $f_0(P)=2p$ and $f_1(P)=|E(G(P))|$. Let $f_{-1}(P)=1$. 

Our main results are the following formulas, which connect the numbers $f_i(P)$ to the numbers  $b_i(J(P))$. 
\begin{thm} \label{main}
Let $P$ be a fixed poset, $p:=|P|$. Let $L:=J(P)$ denote the distributive lattice of order ideals of $P$. Let $k$ be the Sperner number of the poset $P$, i.e., the maximum of the cardinalities of antichains of $P$. Then
\begin{equation} \label{alap}
f_{2p-i-1}(P)=\sum_{m=0}^k (-1)^{m} b_m(L){p-m \choose 2p-i} 
\end{equation}
for each $p\leq i\leq 2p$.
\end{thm}

For example, let $P$ be an antichain with $|P|=p$. Clearly the Sperner number of this poset $P$ is $p$. 
It can be shown easily that $L:=J(P)$ is isomorphic to the Boolean lattice of rank $p$. The vertex set of the graph 
$G(P)$ is the disjoint union $\{x_1,\ldots ,x_p\}\cup \{y_1,\ldots ,y_p\}$, the edge set of $G(P)$ is 
\begin{equation} \label{edge2}
 \{\{x_i,x_j\}:~ 1 \leq i<j\leq p\}\cup  \{\{y_i,y_j\}:~ 1 \leq i<j\leq p\}\cup \{\{x_i,y_j\}:~ i\neq j,1\leq i,j\leq p\}.
\end{equation}
It can be shown easily that $b_i(L)={p \choose i}2^{p-i}$ for each $0\leq i\leq p$ and $f_i(P)={p \choose i+1}2^{i+1}$ for each $-1\leq i\leq p-1$. Hence equation (\ref{alap}) becomes 
$$
{p\choose 2p-i}2^{2p-i}=\sum_{m=0}^p (-1)^m {p \choose m}{p-m\choose 2p-i}2^{p-m}
$$
for each $p\leq i\leq 2p$.

We say that a polytope $Q$ is {\em simplicial}, if all proper faces of $Q$ are simplices.

Let $Q$ be a $d$-dimensional simplicial  polytope. 
We define the $f$-vector of $Q$ as follows:
$$
f(Q):=(f_{-1}(Q),f_0(Q),\ldots ,f_{d-1}(Q))\in {\N}^{d+1},
$$
where
$f_i(Q)$ is the number of $i$--dimensional faces of $Q$ and $f_{-1}(Q)=1$. 
The $h$-vector of $Q$: 
$$
h(Q):=(h_0(Q),\ldots ,h_d(Q))\in {\N}^{d+1},
$$
where
$$
h_k(Q):=\sum_{i=0}^k (-1)^{k-i} {d-i \choose d-k}f_{i-1}(Q)
$$
for each $0\leq k\leq d$.

In particular, $h_0(Q)=1$, $h_1(Q)=f_0(Q)-d$, and  
$$
h_d(Q)=f_{d-1}(Q)-f_{d-2}(Q)\pm \ldots +(-1)^{d-1}f_0(Q)+(-1)^d.
$$

The following equations describe the complete set of linear equations for the coordinates of the $f$-vector of a simplicial polytope.
\begin{thm} \label{Dehn} (Dehn--Sommerville equations)
Let $f(Q)=(f_{-1}(Q),f_0(Q),\ldots ,f_{d-1}(Q))$ be the $f$--vector 
of a $d$--dimensional simplicial polytope. Then 
\begin{equation} \label{DS}
f_{k-1}(Q)=\sum_{i=k}^d (-1)^{d-i} {i \choose k} f_{i-1}(Q).
\end{equation}
for each $0\leq k\leq d$.
\end{thm}

We can write in the following short form these equations:
$$
h_k(Q)=h_{d-k}(Q)
$$
for each $0\leq k\leq d$.

An important special case is the following {\bf Euler--Poincar\'e formula:} 
\begin{equation} \label{EQ}
f_0(Q)-f_{1}(Q)+\ldots +(-1)^d f_{d-1}(Q)=1-(-1)^d.
\end{equation}

The history of the Dehn--Sommerville equations starts 
with M. Dehn, who proved the case $d=5$ in \cite{D}. 
Later Sommerville in \cite{So} proved the general case for simplicial polytopes. V. Klee in \cite{K} gave an elementary proof on the level of simplicial semi--Eulerian complexes. This class includes all the triangulated manifold (without boundary). 

Here we reprove these equations in the special case of simplicial polytopes.   

The outline of the present paper is the following.

First in Chapter 2 we collected the preliminary definitions and results about simplicial complexes, Stanley--Reisner rings, graphs, free resolutions and the Hibi ideal of the poset $P$. In  Chapter 3 we provide 
a short proof for our Theorem \ref{main}. In Chapter 4 we reprove the Dehn--Sommerville equations (\ref{DS}) and give an application using a formula of Peskin and Szpiro \cite{PS}. In our proof we use the homological algebra of free resolutions and the theory of Stanley--Reisner rings. Our results are based on the results of H. Hibi and J. Herzog in \cite{HH}.

\section{Preliminaries}

\subsection{Simplicial complexes}

We say that $\Delta\subseteq 2^{[n]}$ is a {\em simplicial complex} on the vertex set $[n]=\{1,2,\ldots ,n\}$, if  $\Delta$ is a set of subsets of $[n]$ such that $\Delta$ is a down--set, that is, $G\in\Delta$ and $F\subseteq G$ implies that $F\in \Delta$, and $\{i\}\in \Delta$ for all $i$.

The elements of $\Delta$ are called {\em faces} and the {\em dimension} of a face is one less than its cardinality. An $r$-face is an abbreviation for an $r$-dimensional face.
The dimension of $\Delta$ is the dimension of a maximal face.
We use the notation $\mbox{dim}(\Delta)$ for the dimension of $\Delta$.

Let $f_i(\Delta)$ denote the number of $i$--faces of $\Delta$. If $\mbox{dim}(\Delta)=d-1$, then the $(d+1)$--tuple $(f_{-1}(\Delta),\ldots ,f_{d-1}(\Delta))$ is called the {\em $f$-vector} of $\Delta$, where $f_i(\Delta)$ denotes the number of $i$--dimensional faces of $\Delta$. 

For example, let $Q$ be a simplicial polytope. 
The boundary complex $\Delta(Q)$ is formed by 
the set of vertices of all proper faces of $Q$.

A {\em flag} complex is a simplicial complex 
with the property that every
minimal nonface has precisely two elements.

Let $\cF\subseteq 2^{[n]}$ be an arbitrary set system. Define the complement of $\cF$ as
$$
{\cF}':= 2^{[n]}\setminus \cF.
$$
Consider the following set system
$$
\mbox{co}(\cF):=\{[n]\setminus F:~ F\in \cF\}.
$$
We denote by $\cF^*$ the {\em Alexander dual  of $\cF$} 
$$
\cF^*:=\mbox{co}({\cF}')=(\mbox{co}({\cF}))'\subseteq 2^{[n]}.
$$

We collect here some definition from the theory of Stanley--Reisner rings.

Let $\Q$ denote the rational field.  Let $R$ stand for the polynomial ring
$\Q[x_1,\ldots ,x_n]$. We denote by $\Q[x_1,\ldots,x_n]_{\leq s}$ the vector space of all polynomials
over $\Q$ with degree at most $s$.

Let  $\Delta$ be an arbitrary simplicial complex. We associate {\em the Stanley--Reisner ideal} $I(\Delta)$ to the simplicial complex $\Delta$:
$$
I(\Delta):=\langle x_F:~ F\notin \Delta \rangle\unlhd R.
$$ 
Clearly $I(\Delta)$ is a monomial ideal.

The {\em Stanley-Reisner ring} of a simplicial complex $\Delta$ is the quotient ring 
$$
{\Q}[\Delta]:=R/I(\Delta).
$$

Let $I$ be an arbitrary ideal of $R=\Q[x_1,\ldots ,x_n]$. The {\em Hilbert
function} of the algebra $R/I$ is the sequence $h_{R/I}(0), h_{R/I}(1),
 \ldots $. Here $h_{R/I}(m)$ is the dimension over $\Q$ of the factor-space
 $\Q[x_1,\ldots,x_n]_{\leq m}/(I\cap\Q[x_1,\ldots,x_n]_{\leq m})$ (see
\cite[Section 9.3]{BW}).

On the other hand, if we know the $f$-vector of $\Delta$, then we can compute easily the Hilbert function $h_{{\Q}[\Delta]}(t)$ of the Stanley--Reisner ring ${\Q}[\Delta]$.

\begin{lemma} \label{Stan} (Stanley, see Theorem 5.1.7 in \cite{BH}) The Hilbert function of the Stanley--Reisner ring $\Q[\Delta]$ of a $(d-1)$--dimensional simplicial complex $\Delta$ is
\begin{equation} \label{Hilbertfg}
 h_{{\Q}[\Delta]}(t)=\sum_{j=0}^{d-1} f_j(\Delta) {t-1\choose j}.
\end{equation}
\end{lemma} 

Let $\Delta^*$ denote the Alexander dual of the simplicial complex $\Delta$.
We can easily compute $f^*(\Delta)$, the $f$-vector of $\Delta^*$:
\begin{lemma} \label{fvector}
Let $f(\Delta)=(f_{-1}(\Delta),\ldots ,f_{d-1}(\Delta))$ be the $f$--vector of a $(d-1)$--dimensional simplicial complex $\Delta$. Then the $f$--vector of the simplicial complex $\Delta^*$ is: 
\begin{equation}
f^*(\Delta)=f((\Delta)^*)=[\underbrace{1}_{f^*_{-1}},\underbrace{n\choose 1}_{f^*_{0}},\ldots, ,\underbrace{n\choose n-d-1}_{f^*_{n-d-1}},\underbrace{{n\choose n-d}-f_{d-1}}_{f^*_{n-d}},\ldots, \underbrace{{n\choose 2}-f_1}_{f^*_{n-2}}].
\end{equation}
\end{lemma}
\begin{cor} \label{Alex}
Let $f(\Delta)=(f_{-1}(\Delta),\ldots ,f_{d-1}(\Delta))$ be the $f$--vector of a $(d-1)$--dimensional simplicial complex $\Delta$. Then the Hilbert function $h_M(t)$ of the quotient ring $M={\Q}[\Delta^*]=R/I({\Delta}^*)$ is 
\begin{equation} \label{Hfg2}
h_M(t)=h_{{\Q}[\Delta^{*}]}(t)=\sum_{i=0}^{n-d-1} {n\choose i}{t\choose i}+\sum_{j=2}^d \Big({n\choose j}-f_{j-1}\Big) {t\choose n-j}.
\end{equation}
\end{cor}

\proof
If we substitute the $f$--vector $f^*(\Delta)$ of the Alexander dual $\Delta^*$ from Lemma \ref{fvector} to the equation (\ref{Hilbertfg}), we get the result.
\qed

\subsection{Graph theory}

Let $G$ be a finite graph on the vertex set $[n]=\{1,2,\ldots ,n\}$ with no loops and no multiple edges. We will assume in the following that $G$ possesses no isolated vertex. Let $R=\Q[x_1,\ldots ,x_n]$ denote the polynomial ring in 
$n$ variables over the field $\Q$.

We can associate a useful ideal $I(G)$ to the graph $G$.
The {\em edge ideal} of $G$ is the ideal $I(G)$ of $R$ generated by 
the squarefree quadratic monomials $x_ix_j$ such that $\{i,j\}$ is an edge of $G$. 

A finite graph $G$ is {\em bipartitate} if there is a partition $[n]=T\cup T'$ such that each edge of $G$ is of the form $\{j,k\}$, where $j\in T$ and $k\in T'$. It is a well--known fact from graph theory that a finite graph $G$ is bipartitate if and only if $G$ possesses no cycle of odd length.  

The {\em complementary graph} of $G=(V,E)$ is the graph $\overline{G}$ with the vertices of $V$ and edges all the couples $\{v_i,v_j\}$ such that
$i\neq j$ and $\{v_i,v_j\}\notin E$.

A {\em clique} of a graph $G$ is a complete subgraph of $G$. We can associate to a graph $G$ the {\em clique complex} $\Delta(G)$: this  is the collection of all the cliques of the graph $G$, which forms  a simplicial complex.

The following Lemma is an easy consequence of the definitions.
\begin{lemma} \label{graph}
Let $G$ be a simple graph. Then 
\begin{equation}
I(\overline{G})=I(\Delta(G)). 
\end{equation}
\end{lemma}
 
\subsection{Free resolutions}

We introduce some terminology for describing free resolutions.
 
Let $\Q$ denote the rational field.
Let $R$ be the graded ring $\Q[x_1,\ldots ,x_n]$. 
The vector space $R_s=\Q[x_0,\ldots ,x_n]_s$ consists of the homogeneous polynomials of total degree $s$, together with $0$.

Recall that $M$ over $R$ is a {\em graded module} with a family of subgroups 
 $\{M_t:~ t\in {\Z}\}$ of the additive group, where 
$M_t$ are the homogeneous elements of degree $t$, if we can write $M$ in the form
$$
M=\bigoplus_{t\in \Z} M_t
$$
and 
$$
R_sM_t \subseteq M_{s+t}
$$
for all $s\geq 0$ and $t\in \Z$.
If $M$ is finitely generated, then it can be shown easily that $M_t$ are finite dimensional vector spaces over $\Q$. 

Let $M$ be a graded $R$-module and let $d \in \Z$ be an arbitrary integer. We can define
$$
M(d):=\bigoplus_{t \in \Z} M(d)_t,
$$
where $M(d)_t:=M_{d+t}$. Then $M(d)$ is again a graded $R$-module. 

Consider the graded free modules of the form  $R(d_1)\oplus \ldots \oplus R(d_n)$ for any integers $d_1, \ldots , d_n$. We say that these free modules are the {\em twisted graded free modules}.

Let $M$ be a graded $R$--module. A {\em graded resolution} of $M$ is a 
resolution of the form  
\begin{equation}
0\longrightarrow F_n \longrightarrow \ldots \longrightarrow F_1 \longrightarrow M \longrightarrow 0,
\end{equation}
where each $F_l$ is a twisted graded free module  
and each homomorphism $\phi_l:F_l \longrightarrow F_{l-1}$ 
is a graded homomorphism such that $\phi(F_l)_t\subseteq (F_{l-1})_t$ for all  $t\in \Z$.\\

It is a well--known fact from the theory of free resolutions that   every finitely generated $R$--module has a finite graded resolution of length at most $n$  (see \cite[Chapter 6, Theorem 3.8]{CLO2}).

We say that the resolution 
\begin{equation} \label{gradd}
0\longrightarrow F_n \longrightarrow \ldots \longrightarrow F_1 \longrightarrow M \longrightarrow 0
\end{equation}
is {\em minimal} iff 
$\phi_l:F_l \longrightarrow F_{l-1}$ takes the standard basis of 
$F_l$ to a minimal generating set of $im(\phi_l)$ 
for each $l\geq 1$. 

Let $M$ be a finitely generated graded $R$--module.  
Then we define the {\em Hilbert function} $H_M(t)$ by 
$$
H_M(t):=\mbox{dim}_{\Q}\ M_t.
$$

Now we specialize this definition for the case of the homogeneous ideals. 

Let $I\unlhd R$ be a homogeneous ideal of $R$. Then the quotient ring $R/I$ has a natural graded module structure, set $(R/I)_t:=R_t/I_t$, where $I_t:=I\cap R_t$. Thus it comes out from the definitions that if $M:=R/I$ is the quotient graded $R$-module, then $H_M(t)=h_{R/I}(t)$ for each $t\geq 0$.

In the following Theorem we connect the computation of the Hilbert function $H_M(t)$ to the computation of the dimensions of the free graded modules in a graded resolution of $M$.

\begin{thm} \label{Hilbert} (\cite[Chapter 6, Proposition 4.7]{CLO2}) Let $M$ be a graded $R$-module with the graded free resolution   
\begin{equation}
0\longrightarrow F_n \longrightarrow \ldots \longrightarrow F_1 \longrightarrow M \longrightarrow 0.
\end{equation}
If each $F_j$ is the twisted free graded module  $F_j=\bigoplus_{i=1}^{{\beta}_j} R(d_{i,j})$, then 
\begin{equation}
H_M(t)=\sum_{j=1}^k (-1)^j \sum_{i=1}^{{\beta}_j} {n+d_{i,j}+t \choose n}.
\end{equation}
\end{thm}

The numbers ${\beta}_j$ are the {\em Betti numbers} of the module $M$. 

Let $I\subseteq R$ be an arbitrary graded ideal
with graded minimal free resolution
$$
0\longrightarrow \bigoplus_{j=1}^{\beta_1} R(-a_{sj})\longrightarrow \ldots \longrightarrow  \bigoplus_{j=1}^{\beta_s} R(-a_{1j})\longrightarrow R/I \longrightarrow 0
$$
Suppose that 
$$
\mbox{height}(I)=h.
$$
Denote by $e(R/I)$ the Hilbert--Samuel multiplicity of the ring $R/I$. Then by a formula of Peskine and Szpiro \cite{PS}
\begin{equation} \label{PeS}
e(R/I)=\frac{(-1)^i}{h!}\sum_{j=1}^s (-1)^i \sum_{j=1}^{\beta_i} (a_{ij})^h.
\end{equation}

\subsection{Hibi ideals of a poset $P$}

We give here a short summary about the results of H. Hibi and J. Herzog (see \cite{HH}).

Let $P$ be a finite poset, $|P|=p$.
Let $\Q$ denote the rational field. Consider
$$
S:={\Q}[\{x_p,y_p\}_{p\in P}],
$$
the polynomial ring in $2p$ variables.

Let $K\subseteq P$ be an arbitrary order ideal of $P$.
We associate with $K$ the square--free monomial
$$
u_K:=\prod_{p\in K}x_p \prod_{p\in P\setminus K}y_p\in S.
$$
In particular, $u_P:=\prod_{p\in P} x_p$ and $u_{\emptyset}=\prod_{p\in P} y_p$.  

H. Hibi and J. Herzog defined in \cite{HH} the {\em Hibi ideal} 
$$
H(P):=\langle u_K:~ K\in J(P) \rangle\unlhd S
$$ 
of $P$, which is generated by all $u_K$.

They described the following beautiful graded free resolution of $H(P)$ (see Theorem 2.1 of \cite{HH}).

\begin{thm} \label{Hibi}
Let $P$ be an arbitrary poset with $|P|=p$ and denote by $L:=J(P)$ the distributive lattice of order ideals of $P$. 

Let $S:={\Q}[\{x_p,y_p\}_{p\in P}]$
denote the polynomial ring in $2p$ variables. Let $H(P)\unlhd S$ denote the Hibi ideal of $P$. Then $H(P)$ has the following $F_P$ graded minimal free $S$-resolution:
$$
F_P:0\longrightarrow S(-p-k)^{b_k(L)}\longrightarrow S(-p-k+1)^{b_{k-1}(L)}\longrightarrow\ldots
$$
$$
\longrightarrow S(-p-1)^{b_1(L)}\longrightarrow S(-p)^{b_0(L)}
\longrightarrow H(P)\longrightarrow 0,
$$
where $k$ is the Sperner number of $P$, i.e., the maximum of the cardinalities of antichains of $P$. 
\end{thm}

Let $\Gamma_P$ denote the simplicial complex attached to the squarefree monomial ideal $H(P)$, that is, $H(P)=I(\Gamma_P)$.
H. Hibi and J. Herzog described also the Stanley--Reisner ideal of the $(\Gamma_P)^*$ Alexander dual of $\Gamma_P$ (see Lemma 3.1 of \cite{HH}). 
\begin{lemma} \label{Hibi2}
The Stanley--Reisner ideal of the Alexander dual $(\Gamma_P)^*$ is generated by those squarefree monomials $x_iy_j$ such that $p_i\leq p_j$ in $P$.
\end{lemma}

\section{Proof of Theorem \ref{main}}

We follow the following strategy in our proof.

First we compute the Hilbert function $h_M(t)$ of the quotient module $M:=S/H(P)$ from the graded free resolution of $H(P)$. Then we compute this Hilbert function $h_M(t)$ from the theory of Stanley--Reisner rings. These computations yield to a new equation and if we compare the coefficients of ${t \choose i}$ on both side, then the desired equation (\ref{alap}) follows.
 
Let $M:=S/H(P)$ denote the quotient module of the Hibi ideal $H(P)$.
From Theorem \ref{Hilbert} and Theorem \ref{Hibi} we conclude that the Hilbert function $h_M(t)$ of $M$ is 
$$
h_M(t)={t+2p\choose 2p}+\sum_{i=0}^k (-1)^{i+1}b_i(L) {t+2p-(p+i)\choose 2p}=
$$
\begin{equation} \label{Hfg1}
={t+2p\choose 2p}+\sum_{i=0}^k (-1)^{i+1} b_i(L){t+p-i \choose 2p}.
\end{equation} 

\begin{lemma} \label{klikk}
Let $P=\{q_1,\ldots ,q_p\}$ be a finite poset with $|P|=p$. Define the graph $G(P)$ as follows: let the vertex set of $G(P)$ be the disjoint union $\{x_1,\ldots ,x_p\}\cup \{y_1,\ldots ,y_p\}$. Define the edge set of $G(P)$ as
\begin{equation} \label{edge3}
 \{\{x_i,x_j\}:~ 1 \leq i<j\leq p\}\cup  \{\{y_i,y_j\}:~ 1 \leq i<j\leq p\}\cup \{\{x_i,y_j\}:~ q_i\not\leq q_j\}.
\end{equation}
Denote by $\Gamma_P$ the simplicial complex attached to the squarefree monomial ideal $H(P)$, that is, $H(P)=I(\Gamma_P)$.
Then 
$$
(\Gamma_P)^*=\Delta(G(P)).
$$
\end{lemma}

\proof
We write $G_2(P)$ for the bipartitate graph on the vertex set $\{x_1,\ldots ,x_p\}\cup\{y_1,\ldots ,y_p\}$ whose edges are those $\{x_i,y_j\}$ such that $p_i\leq p_j$ in $P$.  It follows from Lemma \ref{Hibi2} that $I((\Gamma_P)^*)=I(G_2(P))$, that is, the Stanley--Reisner ideal of $(\Gamma_P)^*$ is the edge ideal of
$G_2(P)$. Since $G(P)=\overline{G_2(P)}$ by definition, hence 
$$
I((\Gamma_P)^*)=I(G_2(P))=I(\overline{G(P)})=I(\Delta(G(P))),
$$
where we applied Lemma \ref{graph} in the last equality. This means that $(\Gamma_P)^*$ is the clique complex of the graph $G(P)$. \qed

It follows from Lemma \ref{klikk} that $f(P):=(f_{-1}(P),\ldots ,f_{d-1}(P))$ is the $f$--vector of $(\Gamma_P)^*$. 

Let $\Delta$ stand for $(\Gamma_P)^*$, the Alexander dual of $\Gamma_P$.  Then $\Delta^*=((\Gamma_P)^*)^*=\Gamma_P$ and $M=S/H(P)=S/I(\Gamma_P)=\Q[\Gamma_P]={\Q}[\Delta^*]$.

Let $d:=\mbox{dim}((\Gamma_P)^*)+1$. Clearly $n=f_0((\Gamma_P)^*)=2p$. 
It is easy to verify that $d=p$.  Namely
let $K$ denote one of the maximal $d$--clique in the graph $G(P)$. 
Then $d\geq p$ follows from the definition of $G(P)$. On the other hand, if
$V(K)\cap X=\{x_{i_1},\ldots ,x_{i_r}\}$, then $\{y_{i_1},\ldots ,y_{i_r}\}\cap (V(K)\cap Y)=\emptyset$, because $(x_i,y_i)\notin E(G(P))$ for each $1\leq i\leq p$, hence 
$$
d=|V(K)|=|V(K)\cap X|+|V(K)\cap Y|\leq p.
$$

We can apply Corollary \ref{Alex} for $\Delta^*$:
\begin{equation} \label{Hfg3}
h_M(t)=h_{{\Q}[\Delta^*]}(t)=\sum_{i=0}^{n-d-1} {n\choose i}{t\choose i}+\sum_{j=2}^d \Big({n\choose j}-f_{j-1}(P)\Big) {t\choose n-j}.
\end{equation}

Hence the equations (\ref{Hfg1}) and (\ref{Hfg3})
imply that 
$$
{t+2p\choose 2p}+\sum_{m=0}^k (-1)^{m+1} b_m(L){t+p-m \choose 2p}=
$$
$$
=\sum_{i=0}^{n-d-1} {n\choose i}{t\choose i}+\sum_{j=2}^d \Big({n\choose j}-f_{j-1}(P)\Big) {t\choose n-j}.
$$

Since $d=p$ and $n=2p$, we get that
$$
{t+2p\choose 2p}+\sum_{m=0}^k (-1)^{m+1} b_m(L){t+p-m \choose 2p}=
$$
$$
=\sum_{i=0}^{p-1} {2p\choose i}{t\choose i}+\sum_{j=2}^p \Big({2p\choose j}-f_{j-1}(P)\Big) {t\choose 2p-j}.
$$
Using the Vandermonde identities (see \cite{GKP}, 169--170)
$$
{t+2p \choose 2p}=\sum_{i=0}^{2p} {2p\choose i}{t\choose i}
$$
and 
$$
{t+p-m \choose 2p}=\sum_{i=0}^{2p} {p-m \choose 2p-i}{t\choose i}
$$
for each $0\leq m\leq k$, we get that
$$
\sum_{i=0}^{2p} {2p\choose i}{t\choose i}+\sum_{m=0}^k (-1)^{m+1}b_m(L)\Big( \sum_{i=0}^{2p} {p-m \choose 2p-i}{t\choose i}\Big )=
$$
$$
=\sum_{i=0}^{p-1} {2p\choose i}{t\choose i}+\sum_{j=2}^p \Big({2p\choose j}-f_{j-1}(P)\Big) {t\choose 2p-j}.
$$
After simplification we get
$$
\sum_{i=p}^{2p} {2p\choose i}{t\choose i}+\sum_{m=0}^k (-1)^{m+1}b_m(L)\Big( \sum_{i=0}^{2p} {p-m \choose 2p-i}{t\choose i}\Big )=
$$

\begin{equation} \label{jobb}
=\sum_{j=2}^p \Big({2p\choose j}-f_{j-1}(P)\Big) {t\choose 2p-j}.
\end{equation}

Let $p\leq i\leq 2p$ be a fixed index and compare the coefficients of ${t\choose i}$ on both side of equation (\ref{jobb}). Since $\{{t\choose i}:~ i\in \N\}$ is a basis of the vector space $\Q[t]$ over $\Q$, hence these coefficients are the same and equation (\ref{alap}) follows. \qed

\section{The proof of the Dehn--Sommerville equations}

Let $Q$ be a $d$--dimensional simplicial polytope and let $\Delta(Q)$ denote the boundary complex of $Q$.

Let $R$ stand for the polynomial ring
$\Q[x_1,\ldots ,x_n]$.
Here $n:=f_0(Q)$.

We put ${\Delta(Q)}^*$ for the Alexander dual of $\Delta(Q)$. 
Denote by $M:=R/I({\Delta(Q)}^*)$ the Stanley--Reisner ring of ${\Delta(Q)}^*$.
 
First we compute the Hilbert function $h_M(t)$ of $M$
from the following graded free resolution. 

\begin{thm}
 (see \cite{MS}, Example 4.12)
The ideal $I({\Delta(Q)}^*)$ has the following minimal graded free resolution:
$$
{\cF}_Q:0\longrightarrow R(-n)^1\longrightarrow R(1-n)^{f_0(Q)}\longrightarrow\ldots
$$
\begin{equation} \label{free}
\longrightarrow R(d-n-1)^{f_{d-2}(Q)}
\longrightarrow R(d-n)^{f_{d-1}(Q)}\longrightarrow I({\Delta(Q)}^*) \longrightarrow 0.
\end{equation}
\end{thm}

It follows from Theorem \ref{Hilbert} that the Hilbert function of $\Q[{\Delta}(Q)^*]$ is 
$$
h_M(t)=h_{\Q[{\Delta}(Q)^*]}(t)={n+t \choose t}+\sum_{i=-1}^{d-1} (-1)^{d-i}f_i(Q){t+n-n+i+1\choose n}
$$

\begin{equation}  \label{Hfug1}
={n+t \choose t}+\sum_{i=-1}^{d-1} (-1)^{d-i}f_i(Q){t+i+1\choose n}.
\end{equation}

Clearly $\Delta(Q)=((\Delta(Q))^*)^*$. If we apply Corollary \ref{Alex} for the simplicial complex $\Delta:=(\Delta(Q))^*$, then we get
\begin{equation}  \label{Hfug2}
h_M(t)=h_{{\Q}[{\Delta}^*]}(t)=\sum_{i=0}^{n-d-1} {n\choose i}{t\choose i}+\sum_{j=0}^d \Big({n\choose j}-f_{j-1}\Big) {t\choose n-j}.
\end{equation}

Hence the equations (\ref{Hfug1}) and (\ref{Hfug2}) imply that
$$
{n+t \choose t}+\sum_{i=-1}^{d-1} (-1)^{d-i}f_i(Q){t+i+1\choose n}=
$$
$$
=\sum_{i=0}^{n-d-1} {n\choose i}{t\choose i}+\sum_{j=2}^d \Big({n\choose j}-f_{j-1}(Q)\Big) {t\choose n-j}.
$$
Using the Vandermonde identities (see \cite{GKP}, 169--170)
$$
{t+n \choose n}=\sum_{i=0}^{n} {n\choose i}{t\choose i}
$$
and 
$$
{t+i+1 \choose n}=\sum_{j=0}^{i+1} {t \choose n-j}{i+1 \choose j}
$$
for each $i\geq 0$, we get
$$
\sum_{i=0}^{n} {n\choose i}{t\choose i}+\sum_{i=-1}^{d-1} (-1)^{d-i}f_i(Q)(\sum_{j=0}^{i+1} {t \choose n-j}{i+1 \choose j})= 
$$
$$
=\sum_{i=0}^{n-d-1} {n\choose i}{t\choose i}+\sum_{j=0}^d \Big({n\choose j}-f_{j-1}(Q)\Big) {t\choose n-j}.
$$
After simplification we conclude that 
$$
\sum_{j=n-d}^{n} {n\choose j}{t\choose j}+\sum_{i=-1}^{d-1} (-1)^{d-i}f_i(Q)\Big( \sum_{j=0}^{i+1} {t \choose n-j}{i+1 \choose j} \Big)= 
$$
$$
=\sum_{j=0}^d \Big({n\choose j}-f_{j-1}(Q)\Big) {t\choose n-j}.
$$
Consequently
$$
\sum_{j=n-d}^{n} {n\choose j}{t\choose j}+\sum_{j=0}^{d} {t \choose n-j}\Big( \sum_{i=j-1}^{d-1} (-1)^{d-i}{i+1\choose j}f_i(Q) \Big)=
$$
\begin{equation}  \label{Hfug3}
=\sum_{j=0}^d \Big({n\choose j}-f_{j-1}(Q)\Big) {t\choose n-j}.
\end{equation} 

Since 
$$
\sum_{j=n-d}^{n} {n\choose j}{t\choose j}=\sum_{j=0}^{d} {n\choose j}{t\choose n-j},
$$
hence simplifying the equation (\ref{Hfug3}), we get that
$$
\sum_{j=0}^d f_{j-1}(Q){t\choose n-j}+\sum_{j=0}^d {t\choose n-j}\Big(\sum_{i=j-1}^{d-1} (-1)^{d-i} {i+1 \choose j}f_i(Q)\Big)=0.
$$

\begin{equation} \label{veg}
\sum_{j=0}^d {t\choose n-j}\Big(f_{j-1}(Q)+ \sum_{i=j-1}^{d-1} (-1)^{d-i} {i+1 \choose j}f_i(Q)\Big)=0.
\end{equation}

Now we can compare the coefficients of ${t\choose i}$ on both side of equation (\ref{veg}). We can use again the basis property of  $\{{t\choose i}:~ i\in \N\}$. This implies that these coefficients are the same. Thus
$$
f_{j-1}(Q)+ \sum_{i=j-1}^{d-1} (-1)^{d-i} {i+1 \choose j}f_i(Q)=0
$$
for each $0\leq j\leq d$ and the equations (\ref{DS}) follow.

Finally we give an application using a formula of Peskin and Szpiro (see \cite{PS}).

Let $Q$ be a flag simplicial polytope, that is, a simplicial polytope such that the boundary complex of $Q$ is flag.

Let $\Delta(Q)$ denote this boundary 
complex of $Q$. 
Since $Q$ was a flag polytope, the Stanley--Reisner ideal attached to the Alexander dual $\Delta(Q)=((\Delta(Q))^*)^*$ is generated by the monomials $x_px_q$, where $\{p,q\}\notin \Delta(Q)$, hence 
$$
I((\Delta(Q))^*)=\bigcap_{\{p,q\}\notin \Delta(Q)} (x_p,x_q). 
$$
Therefore the squarefree monomial ideal 
$I(\Delta(Q)^*)$ is of height $2$ and the multiplicity
of ${\Q}[\Delta(Q)^*]$ is given by 
$$
e({\Q}[\Delta(Q)^*])=|\{(p,q):~ \{p,q\}\notin \Delta(Q)\}|.
$$
But it is easy to verify that
$$
|\{(p,q):~ \{p,q\}\notin \Delta(Q)\}|={f_0(Q) \choose 2}-f_1(Q).
$$

\begin{cor}
Let $Q$ be a flag simplicial complex. Then
$$
2\cdot \lbrack {f_0(Q) \choose 2}-f_1(Q)\rbrack =\sum_{i=1}^{d+1} (-1)^{i} f_{d-i}(Q)(f_0(Q)-d+i-1)^2
$$
\end{cor}
\proof 
Let $I:=I((\Delta(Q))^*)$ be the Stanley--Reisner ideal attached to the Alexander dual $(\Delta(Q))^*$. Now apply the formula (\ref{PeS}) for the free graded resolution (\ref{free}).

\noindent
{\bf Acknowledgments.} I am indebted to Jonathan Farley for his useful remarks.

\end{document}